\documentclass{amsart}

\numberwithin{equation}{section}

\usepackage{amssymb,amscd}
\usepackage{enumerate, xspace}
\usepackage{hyperref}

\newtheorem{thm}{Theorem}[section]
\newtheorem{lem}[thm]{Lemma}

\newtheorem{prop}[thm]{Proposition}

\newtheorem{rem}[thm]{Remark}

\newcommand\B{{\mathcal B}}
\newcommand\Co{{\mathcal C}}
\newcommand\D{{\mathcal D}}
\newcommand\E{{\mathcal E}}

\newcommand\Lp{{\mathcal L}}

\newcommand\M{{\mathcal M}}

\newcommand\boN{{\mathcal{N}}}
\newcommand\LL[2]{\Lambda(#1,#2)}

\newcommand\C{{\mathbb C}}

\newcommand\N{{\mathbb N}}

\newcommand\R{{\mathbb R}}

\newcommand\Z{{\mathbb Z}}

\newcommand\ve{\varepsilon}
\newcommand\eps{\varepsilon}
\newcommand\vf{\varphi}
\renewcommand{\epsilon}{\ve}

\newcommand{\p}{\mathbf{p}}
\newcommand{\q}{\mathbf{q}}
\newcommand{\x}{\mathbf{x}}
\newcommand{\y}{\mathbf{y}}
\newcommand{\ii}{\mathbf{i}}
\newcommand{\jj}{\mathbf{j}}
\newcommand{\kk}{\mathbf{k}}

\DeclareMathOperator{\Var}{Var} \DeclareMathOperator{\var}{var}
\DeclareMathOperator{\Id}{Id} \DeclareMathOperator{\dist}{dist}
 \DeclareMathOperator{\Lip}{Lip}
\newcommand{\moins}{\backslash}
\newcommand{\dd}{\, {\rm d}}
\newcommand{\norm}[1]{\left\| #1 \right\|}
\newcommand{\st}{\,:\,}

\newcommand{\nsigma}{\alpha}

\makeatletter
\newcommand{\subjclassname@NEW}{2000 Mathematics Subject Classification}
\makeatother

\begin{document}
\title{Limit theorems for coupled interval maps} \author{Jean-Baptiste Bardet,
  S\'{e}bastien Gou\"{e}zel and Gerhard Keller} \address{J.-B. Bardet and S.
  Gou\"{e}zel: IRMAR/UFR Math\'{e}matiques, Universit\'{e} de Rennes 1, Campus de
  Beaulieu, 35042 Rennes Cedex, France;\; G. Keller: Mathematisches Institut,
  Universit\"{a}t Erlangen-N\"{u}rnberg, Bismarckstr. 1 1/2, 91054 Erlangen,
  Germany} \email{ {\tt
    jean-baptiste.bardet@univ-rennes1.fr,sebastien.gouezel@univ-rennes1.fr,
    keller@mi.uni-erlangen.de}} \date{\today} \thanks{G.K. thanks the
  colleagues at the UFR Math\'{e}matiques of the University of Rennes 1 for
  their hospitality during his stay in March and April 2006.}
\begin{abstract}
  We prove a local limit theorem for Lipschitz continuous observables on a
  weakly coupled lattice of piecewise expanding interval maps. The core of the
  paper is a proof that the spectral radii of the Fourier-transfer operators
  for such a system are strictly less than $1$. This extends the approach of
  \cite{kl:coupling2} where the ordinary transfer operator was studied.
\end{abstract}
\keywords{Coupled map lattice, piecewise expanding map, spectral
gap, local limit theorem} \subjclass[NEW]{37L60,60F05}

\maketitle
\section{Results}

This paper deals with the issue of probabilistic limit theorems in
dynamical systems, i.e., limit theorems for the Birkhoff sums $S_n
f=\sum_{k=0}^{n-1} f\circ T^k$, where $T$ is a probability
preserving transformation of a space $X$ and $f:X\to \R$ is an
appropriate measurable function. There are currently many techniques
available to prove the central limit theorem $S_n f/\sqrt{n} \to
\boN(0,\sigma^2)$, let us mention for example elementary techniques,
martingales, spectral arguments. On the other hand, if one is
interested in the local limit theorem $\mu\{ S_n f\in [a,b]\} \sim
\frac{|b-a|}{\sigma \sqrt{2\pi n}}$, the scope of possible
techniques is much more narrow: all known proofs rely on spectral
analysis of transfer operators. Therefore, the class of systems for
which a local limit theorem is proved is much smaller.

We are interested in limit theorems for coupled map lattices. The
only previous result in this context is \cite{barde1}, where central
limit theorem, moderate deviations principle and a partial large
deviations principle were established under strong analyticity
assumptions on the local map and the coupling.  In this paper, we
establish central and local limit theorems for coupled interval maps
under much weaker assumptions. More precisely, we study the same
class of systems as in \cite{kl:coupling2}. We emphasize on local
limit theorem, since it is the most demanding result. But our
method, relying on spectral analysis of transfer operators, gives
other limit theorems, see Remark \ref{rem:moreres} below.

Let us recall the setup from \cite{kl:coupling2}. Given a compact
interval $I\subset \R$ we will consider the phase space
$\Omega:=I^{\Z^d}$.  In the following we always assume without loss
of generality that $I=[0,1]$.

The single site dynamics is given by a map $\tau: I\to I$. We assume
$\tau$ to be a continuous, piecewise $\Co^2$ map from $I$ to $I$
with singularities at $\zeta_1,\dots,\zeta_{N-1}\in(0,1)$ in the
sense that $\tau$ is monotone and $\Co^2$ on each component of
${I}\setminus\{\zeta_0=0,\zeta_1,\ldots,\zeta_{N-1},\zeta_N=1\}$.
We assume that $\tau', \tau''$ are bounded and that $\inf|\tau'|>2$.
Next, we define the unperturbed dynamics $T_0:\Omega\to\Omega$ by
$[T_0(x)]_{\p}:=\tau(x_\p)$.

To define the perturbed dynamics we introduce couplings
$\Phi_\eps:\Omega\to\Omega$ of the form
$\Phi_\eps(\x):=\x+A_\eps(\x)$. We say that such a coupling has
range $r$ and strength $\eps$ if for all $\kk,\p,\q\in {\Z^d}$
\begin{equation}
  \label{eq:couplingass}
  |(A_\eps)_\p|_\infty\le 2\eps,\quad
  |(DA_\eps)_{\q\p}|_\infty\le 2\eps,\quad
  |\partial_{\kk}(DA_\eps)_{\q\p}|_\infty\le 2\eps\;,
\end{equation}
and $\partial_\p\Phi_{\eps,\q}=0$ whenever $|\p-\q|>r$. The
diffusive nearest neighbor coupling used in \cite{MH-93}, and in
much of the numerical literature, is defined by
\begin{equation}
  \label{eq:diffusive-nnc}
  [\Phi_\eps (\x)]_\p
  =
  x_{\p}+\frac{\eps}{2d}\sum_{|\p-\q|=1} (x_\q-x_\p) \quad(\p\in
  {\Z^d})\;,
\end{equation}
and it is a trivial example of such a coupling with range $r=1$ and
strength $\eps$. The dynamics $T_\eps:\Omega\to\Omega$ that we wish
to investigate is then defined as
  \begin{equation}
  T_\eps:=\Phi_\eps\circ T_0\;.
  \end{equation}
Let $m$ denote Lebesgue measure on the interval $I$. The following
result is proved in \cite{kl:coupling2}:
\begin{thm}
\label{Thm:InvMeas} For each $r\in \N$, there exists $\eps_0(r)>0$
such that, for any coupling $\Phi_\eps$ of range $r$ and strength
$0\leq\eps\leq\eps_0(r)$, there exists a unique measure $\mu_\eps$
such that, for $m^{\otimes \Z^d}$-almost every point $\x$,
  \begin{equation}
  \label{eq:ConvBirkhoof}
  \frac{1}{n}\sum_{k=0}^{n-1} \delta_{T_\eps^k \x} \to \mu_\eps.
  \end{equation}
  This measure $\mu_\eps$ has in fact many additional properties: it is the
  unique invariant measure in the class $\B$ of measures of bounded variation
  (that we will define later), it is exponentially mixing both in time and
  space, and the convergence \eqref{eq:ConvBirkhoof} holds for $\mu$-almost
  every point whenever $\mu$ is a measure of bounded variation.
\end{thm}

In this paper, we prove the following theorem.
\begin{thm}
\label{MainThm} For each $r\in \N$, there exists $0<\eps_1(r)\leq
\eps_0(r)$ satisfying the following property. Let $\Phi_\eps$ be a
coupling of range $r$ and strength $0\leq\eps\leq\eps_1(r)$, and let
$\mu_\eps$ denote the corresponding invariant measure given by
Theorem \ref{Thm:InvMeas}. Let $f: \Omega \to \R$ be a Lipschitz
function depending on a finite number of coordinates, with $\int
f\dd\mu_\eps=0$.
\\
\textbf{Central limit theorem.} There exists $\sigma^2\geq 0$ such
that $\frac{1}{\sqrt{n}} \sum_{k=0}^{n-1} f\circ T_\eps^k$ converges
in distribution to $\boN(0,\sigma^2)$, with respect to the measure
$\mu_\eps$. Moreover, $\sigma^2=0$ if and only if there exists a
measurable function $u:\Omega\to \R$ such that $f=u-u\circ T_\eps$
$\mu_\eps$-almost everywhere.
\\
\textbf{Local limit theorem.} Assume additionally that, whenever
$u:\Omega\to \R$ is measurable and $\lambda\in \R^*$, the function
$f-u+u\circ T_\ve \mod \lambda$ is not $\mu_\eps$-almost everywhere
constant -- we say that $f$ is \emph{aperiodic}. In particular, the
variance $\sigma^2$ in the central limit theorem is nonzero. Then,
for any compact interval $I\subset \R$,
  \begin{equation}
  \sigma\sqrt{2\pi n}\cdot \mu_\eps\{ x \st S_n f(x)\in I\} \to |I|.
  \end{equation}
Here, $|I|$ denotes the length of the interval $I$.
\end{thm}

It is probably possible to weaken the assumptions, by replacing the
finite range interaction by a short range interaction, and by
allowing the function $f$ to depend on all coordinates but with an
exponentially small influence of far away coordinates (by mimicking
the techniques of \cite[Section 5]{kl:coupling2}). On the other
hand, it is unclear whether it is possible to remove the continuity
assumption on $\tau$ (notice that, for finite range interactions,
this condition is not required in \cite{kl:coupling2}).

On the technical level, Theorem \ref{MainThm} is a consequence of a
spectral description of perturbed transfer operators acting on a
suitable Banach space, that we now describe. Denote by $\M$ the set
of complex Borel measures on $\Omega$ where $\Omega$ is equipped
with the product topology.

Let $\Co$ be a set of objects ``acting on functions depending on
finitely many coordinates'', defined as follows. An element of $\Co$
is a family $(\mu_\Lambda)$, where $\Lambda$ goes through the finite
subsets of $\Z^d$, such that $\mu_\Lambda$ is a complex measure on
$I^\Lambda$, and such that if $\Lambda'\subset \Lambda$ then the
projection of $\mu_\Lambda$ on $I^{\Lambda'}$ is $\mu_{\Lambda'}$.
Formally, $\Co$ is the projective limit of the spaces of complex
measures on $I^\Lambda$, $\Lambda$ finite subsets of $\Z^d$. This is
a complex vector space, and we will not use any topology on it. Note
that there is a canonical inclusion of $\M$ in $\Co$. If $u$ is a
bounded measurable function depending on a finite number of
coordinates, and $\mu\in \Co$, then it is possible to define
canonically $u\mu\in \Co$.

If $\mu\in \Co$ and $\vf$ is a bounded measurable function depending
on a finite number of coordinates, it is possible to define
$\mu(\vf)$ as $\mu_\Lambda(\vf)$ whenever $\Lambda$ is large enough.
If $\vf$ depends on finitely many coordinates, then $\vf\circ
T_\eps$ also depends on finitely many coordinates. This implies
that, for any $\mu\in \Co$, there exists a unique $\nu\in \Co$ such
that, for any $\vf$,
  \begin{equation}
  \nu(\vf)=\mu(\vf\circ T_\eps).
  \end{equation}
We write $\nu=P_\eps \mu$. Thus, $P_\eps$ is a linear operator from
$\Co$ to $\Co$. It is the so-called \emph{transfer operator} of the
map $T_\eps$. The image under $P_\eps$ of a measure is still a
measure.

\begin{thm}
\label{thm:SpectralDescription} There exists a subspace $\D$ of
$\Co$ endowed with a complete norm $\norm{\cdot}$ with the following
properties. First, $\D$ contains the set of measures with bounded
variation, and for any $\mu\in \D$, $|\mu(1)| \leq \norm{\mu}$.
Moreover, for any finite subset $\Lambda$ of $\Z^d$, there exists a
constant $C(\Lambda)$ such that, for any $u:\Omega\to \C$ depending
only on coordinates in $\Lambda$ and Lipschitz, for any $\mu\in \D$,
$u\mu$ also belongs to $\D$ and
  \begin{equation}
  \label{eq:Algebra}
  \norm{ u\mu} \leq C(\Lambda) (\Lip(u)+|u|_\infty) \norm{\mu},
  \end{equation}
where $\Lip(u)$ denotes the best Lipschitz constant of $u$. (The
norm $\|\,.\,\|$ is defined in equation \eqref{eq:TotalMass}.)

For any $r\in \N$, there exists $\epsilon_1(r)>0$ such that, if
$\Phi_\eps$ is a coupling with range $r$ and strength
$0\leq\eps\leq\epsilon_1(r)$, then the following holds.
\begin{itemize}
\item
If $\mu \in \D$, then $P_\eps \mu\in \D$ and $\norm{P_\eps \mu}\leq
C\norm{\mu}$ for some constant $C$. In fact, the operator $P_\eps$
has a simple eigenvalue at $1$ and the rest of its spectrum is
contained in a disk of radius $<1$.
\item
Let $f$ be a Lipschitz function depending on a finite number of
coordinates. Then the map $t\mapsto P_{t,\eps} = P_\eps(e^{itf}.)$
is an analytic map from $\R$ to $\Lp(\D)$, the set of continuous
linear operators on $\D$. If $f$ is aperiodic, then the spectral
radius of $P_{t,\eps}$ is $<1$ for any $t\not=0$.
\end{itemize}
\end{thm}

The derivation of Theorem \ref{MainThm} from Theorem
\ref{thm:SpectralDescription} is classical. Note however that some objects
from $\D$ occurring in the proof are not known to be measures, so that one
cannot directly cite \cite{parry-pollicott}, for example.  So we will sketch
the details of the proof in Appendix \ref{append}, because this seems clearer
than applying an abstract result like 
e.g.~\cite[Theorem VII.1.8]{kato:pe} or~\cite[Corollary
III.11]{hennion_herve}.

\begin{rem}
\label{rem:moreres} Theorem \ref{thm:SpectralDescription} implies
even more precise results: the limit theorems of Theorem
\ref{MainThm} hold not only for $\mu_\eps$, but also for any
probability measure $\mu$ which belongs to $\D$ (and in particular
for any probability measure of bounded variation). Additionally,
further refinements of the central limit theorem hold. For example,
the speed of convergence in the central limit theorem is
$O(1/\sqrt{n})$, a renewal theorem holds, as well as a large
deviation inequality (see again \cite{hennion_herve} for further
details). One can also derive in the same way as in \cite{barde1}
the moderate deviations principle.
\end{rem}

The rest of the paper will be devoted to the proof of Theorem
\ref{thm:SpectralDescription}. The main problem will be to get a
Lasota-Yorke inequality with compactness, since the space of
measures of bounded variation is not compact in the space of finite
measures. We will therefore use artificial extensions as in
\cite{kl:coupling2}, but we will lose control in the ``central box''
due to the factor $e^{itf}$. This loss will be compensated by the
fact that, in large but finite boxes, the measures of bounded
variation form a compact subset of the set finite measures.
Technically, we will have to take larger and larger boxes as $t$
increases, but this causes no harm.

\section{Functional analytic constructions}

\subsection{Abstract tools}

We will need the following lemma.
\begin{lem}
\label{lem:StableCohom} Let $T$ be a transformation preserving a
probability measure $\mu$. Let $n>0$. Then a function $f$ is
aperiodic for $T$ if and only if $S_n f$ is aperiodic for $T^n$.
\end{lem}
\begin{proof}
If $f$ is periodic, then there exist $c,d>0$ and $u$ measurable such
that $f=u-u\circ T+d \mod c$. Therefore, $S_n f= u-u\circ T^n + nd
\mod c$, hence $S_nf$ is periodic.

Conversely, assume that $S_n f$ is periodic for $T^n$, i.e., $S_n
f=u-u\circ T^n + d \mod c$. Then $S_n (f-u+u\circ T) = d \mod c$.
For any function $v$, $S_n v$ is cohomologous to $nv$ (since $v\circ
T^k$ is cohomologous to $v$). Therefore, there exists a function $w$
such that
  \begin{equation}
  n(f-u+u\circ T) = S_n (f-u+u\circ T) + w-w\circ T
  = d+ w-w\circ T \mod c.
  \end{equation}
Therefore, $f$ is cohomologous to a constant modulo $c/n$, and $f$
is periodic.
\end{proof}

We will also need the following formula on the essential spectral
radius.
\begin{lem}
\label{lem:hennion} Let $Q$ be a continuous linear operator on a
complex Banach space $(B, \norm{\cdot})$. Assume that there exists a
semi-norm $\norm{\cdot}_w$ on $B$ such that any sequence $x_n$ in
$B$ with $\norm{x_n} \leq 1$ contains a Cauchy subsequence for
$\norm{ \cdot}_w$. Assume moreover that there exist $\sigma>0$ and
$C>0$ such that, for any $x\in B$,
  \begin{equation}
  \norm{Qx} \leq \sigma \norm{x} +C \norm{x}_w.
  \end{equation}
Then the essential spectral radius of $Q$ is at most $\sigma$.
\end{lem}
This is a version of a theorem by Hennion \cite{hennion}, where one
does not need to be able to iterate the operator for the weak norm
(in the forthcoming application, the operator $Q$ will indeed not be
continuous for the weak norm).
\begin{proof}
Let $M>0$ be such that $\norm{Qx}\leq M\norm{x}$. Notice also that
there exists by assumption a constant $C>0$ such that
$\norm{x}_w\leq C \norm{x}$ for all $x\in B$. It allows to define a
new seminorm on $B$ by $\norm{x}_w' = \sum_{n\geq 0} (2M)^{-n}
\norm{ Q^n x}_w$. It satisfies the same compactness assumptions as
$\norm{\cdot}_w$. Moreover, $Q$ is continuous for this seminorm. We
can therefore iterate the equation
  $
   \norm{Qx} \leq \sigma \norm{x} +C \norm{x}'_w,
  $
and get an estimate
  \begin{equation}
  \norm{Q^n x} \leq \sigma^n \norm{x} + C_n \norm{x}'_w.
  \end{equation}
The aforementioned theorem of Hennion \cite[Corollaire 1]{hennion}
gives the conclusion.
\end{proof}

\subsection{Measures of bounded variation}
\label{subsec:BV} The concept of measures of bounded variation will
play a central role. For $\mu\in \M$ and $\p\in \Z^d$, we define
  \begin{equation}
    \Var_\p \mu:=\sup_{|\vf|_{\Co^0(\Omega)}\leq
  1}\mu(\partial_\p\vf)\;.
  \end{equation}
Here, the sup is restricted to functions which are $C^1$ in $x_\p$,
depending only on a finite number of coordinates. Let also
  \begin{equation}\label{eq:variation}
  \Var \mu=\sup_{\p\in \Z^d} \Var_\p \mu\;.
  \end{equation}
The set $\B:=\{\mu\in\M\;:\;\Var\mu<\infty\}$  consists of measures
whose finite dimensional marginals are absolutely continuous with
respect to Lebesgue and the density is a function of bounded
variation. In fact, ``$\Var$'' is a norm and, with this norm, $\B$
is a Banach space.

We define also in the same way, for any subset $\Lambda$ of $\Z^d$
and any measure $\mu_\Lambda$ on $I^\Lambda$,
 \begin{equation}\label{eq:variationsubset}
  \Var_\Lambda \mu_\Lambda=\sup_{\p\in \Lambda} \Var_\p \mu_\Lambda\;.
  \end{equation}

We also need the usual total variation norm on complex measures:
\begin{equation}\label{eq:regularnorm}
|\mu|:=\sup_{|\vf|_{\Co^0(\Omega)}\leq 1}\mu(\vf)\;.
\end{equation}
Just like in \cite[Sect.~3.3]{kl:coupling1} one checks easily that
\begin{equation}
  \label{eq:two-norms}
  |\mu|\leq\frac12\Var_\p\mu\quad(\p\in\Z^d).
\end{equation}

For $\mu\in \M$, let $A(\mu)$ denote its absolute value, it is a
positive measure.
\begin{lem}
\label{lem:AinB} If $\mu\in \B$, then $A(\mu)\in \B$ and $\Var
A(\mu) \leq \Var(\mu)$.
\end{lem}
\begin{proof}
When $\mu$ is a measure with bounded variation on an interval, then
the formula $\Var A(\mu) \leq \Var(\mu)$ is a direct consequence of
the formula
  \begin{equation}
  \Var(\mu) = \inf_{ f\dd m = \dd\mu} \sup_{x_1 <\dots<x_k} \sum
  | f(x_{i+1}) -f(x_i)|\;.
  \end{equation}
Indeed, if $\dd\mu=f\dd m$ then $\dd A(\mu)= |f|\dd m$, and the
formula $| |f|(x_{i+1}) - |f|(x_i)| \leq |f(x_{i+1})-f(x_i)|$
implies the conclusion.

In dimension $n$, the variation of a measure can be written as the
integral of one-dimensional variations (see e.g.~(43) in
\cite{kl:coupling1}). Hence, the result is implied by the
one-dimensional result.

Consider now a measure $\mu\in \B$. If $A(\mu)=0$, there is nothing
to do. Otherwise, we can assume without loss of generality that
$A(\mu)$ is a probability measure. There exists a measurable
function $\vf$, of absolute value almost everywhere equal to one,
such that $\mu=\vf A(\mu)$. Let $\psi$ be a $C^1$ test function
depending on a finite number of coordinates and bounded by $1$, and
let $\q\in \Z^d$. For any finite box $\Lambda$ (containing all the
coordinates on which $\psi$ depends), the finite dimensional result
implies
  \begin{equation}
  \label{eq:ALambda}
  A( \pi_\Lambda \mu) ( \partial_\q \psi) \leq \Var( A(\pi_\Lambda
  \mu)) \leq \Var( \pi_\Lambda \mu) \leq \Var(\mu).
  \end{equation}
Let $\vf_\Lambda$ denote the conditional expectation (for the
measure $A(\mu)$) of the function $\vf$ with respect to the
$\sigma$-algebra of sets depending only on coordinates in $\Lambda$.
Then
  $  \pi_\Lambda(\mu) = \pi_\Lambda( \vf A \mu) = \vf_\Lambda
  \pi_\Lambda(A \mu ). $
Therefore,
  $
  A( \pi_\Lambda \mu ) = | \vf_\Lambda| \pi_\Lambda(A \mu).
  $
Hence, \eqref{eq:ALambda} reads
  \begin{equation}
  \label{partial_conv}
  \int | \vf_\Lambda| \partial_\q \psi \dd (A \mu) \leq \Var(\mu).
  \end{equation}
When the box $\Lambda$ increases, the sequence of functions
$\vf_\Lambda$ converges in $L^1(A \mu)$ to $\vf$, by the martingale
convergence theorem. Therefore, $|\vf_\Lambda|$ converges to
$|\vf|=1$. Taking the limit in \eqref{partial_conv}, we get
  \begin{equation}
  \int \partial_\q\psi \dd (A \mu) \leq \Var(\mu). \qedhere
  \end{equation}
\end{proof}

An element $\mu$ of $\B$ gives canonically rise to an element
$(\mu_\Lambda)_\Lambda$ of $\Co$ by taking the induced measure on
every finite subset of $\Z^d$. It satisfies
  \begin{equation}
  \label{FiniteVar}
  \sup_\Lambda \Var_\Lambda(\mu_\Lambda)<\infty.
  \end{equation}
\begin{lem}
\label{lem:BoundedVarOK} Conversely, consider an element
$(\mu_\Lambda)_\Lambda$ of $\Co$ satisfying \eqref{FiniteVar}. Then
it comes from an element of $\B$.
\end{lem}
\begin{proof}
Let $\Lambda_n$ be an increasing sequence of boxes. Define a measure
$\mu_n$ on $I^{\Z^d}$ by $\mu_n= \mu_{\Lambda_n} \otimes m^{\otimes
\Z^d \moins \Lambda_n}$. The sequence $\mu_n$ has uniformly bounded
variation. Let $\mu$ be one of its weak limits. Its marginal on each
box $\Lambda$ coincides with $\mu_\Lambda$ by construction.
\end{proof}

For $u :\Omega\to \R$ and $\p\in\Z^d$, let
  \begin{equation}
  \Lip_\p(u)= \sup_{ \x \in I^{\Z^d \moins \{\p\}}} \sup_{x_\p \not=
  x'_\p \in I} \frac{ u( x_\p, \x)-u(x'_\p,\x)}{|x_\p - x'_\p|}.
  \end{equation}
\begin{lem}
\label{lem:LipschitzControl} For any $u:\Omega \to \R$ depending on
a finite number of coordinates, any $\mu\in \B$ and any $\p\in
\Z^d$,
  \begin{equation}
  \Var_\p( u\mu) \leq \sup |u| \Var_\p(\mu) + \Lip_\p(u) |\mu|.
  \end{equation}
\end{lem}
\begin{proof}
In one dimension, this is a consequence of \cite[Lemma
2.2(b)]{kl:coupling1} and the fact that a Lipschitz function is
differentiable almost everywhere and is equal to the integral of its
derivative. This extends to finite boxes by (43) in
\cite{kl:coupling1}. Taking the supremum over finite boxes yields
the conclusion of the lemma.
\end{proof}

\subsection{A family of extensions}

For $\p \in \Z^d$, denote by $\B_\p$ the set of measures $\mu$ in
$\B$ such that, whenever a test function $\vf$ does not depend on
the coordinate $\p$, then $ \mu(\vf)=0$.

We can now define a family of extensions. We adapt the construction
of \cite[Section 3]{kl:coupling2}, the main difference being that we
keep a central part of the measure on a finite subset of $\Z^d$.

Let $\Lambda$ be a finite subset of $\Z^d$, we define a space
$\E(\Lambda)$ as follows. An element of $\E(\Lambda)$ is a family
$\mu=(\mu_c, (\mu_{\p})_{p\in \Z^d\moins \Lambda})$ such that
$\mu_c$ is a measure of the form $\nu\otimes m^{\Z^d\moins \Lambda}$
where $\nu$ is a measure on $I^\Lambda$, and $\mu_\p\in \B_\p$.
Here, $m$ denotes Lebesgue measure on $I$. We assume moreover
  \begin{equation}
  \norm{\mu}:= \max( \Var(\mu_c), \sup_{\p \in \Z^d\moins\Lambda}
  \Var(\mu_\p)) <\infty.
  \end{equation}
On $\E(\Lambda)$, we also define a ``weak norm'' by
  \begin{equation}
  \norm{\mu}_w= |\mu_c|.
  \end{equation}
The unit ball of $(\E(\Lambda), \norm{\cdot})$ is relatively compact
for the seminorm $\norm{\cdot}_w$.

There is a canonical projection from $\E(\Lambda)$ to $\Co$, given
by the sum of the measures $\mu_c$ and $(\mu_\p)_{p\in \Z^d\moins
\Lambda}$. We will denote it by $\pi_{\E(\Lambda)}$ or simply by
$\pi$.

We describe now a (non-canonical) redistribution process introduced
in \cite{kl:coupling2}. Let $B$ be a subset of $\Z^d$, of
cardinality $J\in [0,\infty]$. Let $\sigma: [0,J) \to B$ be an
enumeration of the points in $B$. For $j\leq J$, let
$B_j=\sigma[0,j)$. In particular, $B_0=\emptyset$ and $B_J=B$. If
$\mu\in \B$, define measures $\mu_\p$, for $\p\in B$, by
  \begin{equation}
  \mu_\p=\pi_{\Z^d\moins B_j} \mu \otimes m^{\otimes B_j}
  - \pi_{\Z^d\moins B_{j+1}} \mu \otimes m^{\otimes B_{j+1}},\quad
  \text{where }j=\sigma^{-1}(\p).
  \end{equation}
By construction, $\mu=\pi_{\Z^d\moins B}\mu\otimes m^{\otimes B}+
\sum_{\p\in B} \mu_\p$, and $\mu_\p\in \B_\p$ satisfies
$\Var(\mu_\p) \leq 2 \Var \mu$. We say that $\pi_{\Z^d\moins
B}\mu\otimes m^{\otimes B}$ is the part of $\mu$ remaining at the
end of the redistribution process.

Using this process for $B=\Z^d\moins \Lambda$, we obtain a map
$H_\Lambda$ which associates to any $\mu\in \B$ an element
$H_\Lambda(\mu) \in \E(\Lambda)$. It satisfies
$\norm{H_\Lambda(\mu)} \leq 2\Var(\mu)$, and $\pi\circ H_\Lambda=
\Id$.

\smallskip

Finally, let $f$ be a Lipschitz function depending on a finite
number of coordinates, and let $t\in \R$. Assume that the function
$tf$ depends only on coordinates in $\Lambda$. For each $n\in \N$,
we define on $\E(\Lambda)$ an operator $Q_{t,\ve,n,\Lambda}$, which
is a (non-canonical) lift of $P_{t,\ve}^n$ on $\Co$. Starting from
$\mu=(\mu_c, (\mu_\p))\in \E(\Lambda)$, apply first $P_{t,\ve}^n$ to
each measure $\mu_c$ and
$\mu_{\p}$. Then, redistribute the mass as follows:\\
\noindent $\bullet$  For $\dist(\p,\Lambda)>nr$, distribute
$P_{t,\ve}^n\mu_\p$ to $B=\{ \q \st |\q-\p|\leq nr\}$. The points of
$B$ are all outside of $\Lambda$. Moreover, since $\mu_\p\in \B_\p$
and $tf$ depends only on coordinates in $\Lambda$, we have
$\pi_{\Z^d\moins B}(P_{t,\ve}^n\mu_\p)=0$, i.e., there is no mass
remaining at the end of this redistribution process.
\\
\noindent $\bullet$ For the other measures, use
$H_\Lambda$.\\
We get as in \cite[Lemma 3.1]{kl:coupling2}
  \begin{equation}
  \label{ExpandsWeak}
  \norm{ Q_{t,\ve,n,\Lambda} \mu} \leq 2 B(\Lambda,n,r)\sup\left( \Var(P_{t,\ve}^n \mu_c),
  \sup_{\p \in \Z^d\moins \Lambda} \Var(P_{t,\ve}^n \mu_\p)\right)
  \end{equation}
with $B(\Lambda,n,r)=\#\{ \q \in \Z^d \st
  \dist(\q,\Lambda) \leq nr\} + \#\{\q\in \Z^d \st |\q|\leq
  nr\}$, since every new measure receives a
contribution from a number of sites bounded by $B(\Lambda,n,r)$.
Note that we have written $n$ as an index and not an exponent, in
$Q_{t,\ve,n,\Lambda}$, since these operators are not the powers of a
single operator due to the (non-canonical) redistribution process.

Note that the extension $\E(\emptyset)$ is at the heart of the proof
of \cite{kl:coupling2}.

\subsection{Construction of a canonical extension}
\label{subsec:DescribesD}

Let $\Lambda$ be a finite subset of $\Z^d$. Let $\E(\Lambda)_0
\subset \E(\Lambda)$ be the kernel of $\pi_{\E(\Lambda)}$, i.e., the
elements of $\E(\Lambda)$ which induce the zero measure on the
basis. This is a closed subspace of $\E(\Lambda)$, we can therefore
consider the quotient space $\D(\Lambda):= \E(\Lambda)/
\E(\Lambda)_0$ with its canonical norm. The map $\pi_{\E(\Lambda)}$
induces a map $\pi_{\D(\Lambda)}: \D(\Lambda) \to \Co$, which is
injective. In this way, we can therefore consider $\D(\Lambda)$ as a
subspace of $\Co$.

Since $\pi_{\E(\Lambda)}\circ Q_{t,\ve,n,\Lambda}=P_{t,\ve}^n\circ
\pi_{\E(\Lambda)}$, the operator $Q_{t,\ve,n,\Lambda}$ leaves
$\E(\Lambda)_0$ invariant, and induces therefore a map $\bar
Q_{t,\ve,n,\Lambda}$ on $\D(\Lambda)$. An interesting consequence of
this construction is that $\bar Q_{t,\ve,n,\Lambda}= \bar
Q_{t,\ve,1,\Lambda}^n$, i.e., we are really dealing with the powers
of a single operator. This is due to the fact that the
non-canonicity in the redistribution process is killed by the
quotient, any redistribution would induce the same map on
$\D(\Lambda)$.

\begin{prop}
\label{prop:EquivNorms} If $\Lambda, \Lambda'$ are two finite
subsets of $\Z^d$, then the subsets $\pi_{\D(\Lambda)}(\D(\Lambda))$
and $\pi_{\D(\Lambda')}(\D(\Lambda'))$ of $\Co$ are equal, and the
induced norms are equivalent.
\end{prop}
\begin{proof}
It is sufficient to prove this for $\Lambda'\subset \Lambda$.
Consider $\Lambda'\subset \Lambda$, and construct a continuous
linear map from $\E(\Lambda)$ to $\E(\Lambda')$ by redistributing
the mass of $\mu_c$ in any convenient way. This induces a map from
$\D(\Lambda)$ to $\D(\Lambda')$. Conversely, starting from an
element of $\E(\Lambda')$, we can consider $\mu_c+\sum_{\p \in
\Lambda \moins \Lambda'}\mu_\p$ and redistribute it in any way, to
get an element of $\E(\Lambda)$. Going to the quotient gives a
canonical map from $\D(\Lambda')$ to $\D(\Lambda)$, which is inverse
to the previous one. Hence, we have constructed a canonical
isomorphism between $\D(\Lambda)$ and $\D(\Lambda')$, which commutes
with the projections $\pi_{\E(\Lambda)}$ and $\pi_{\E(\Lambda')}$.
We get the proposition by projecting everything in $\Co$.
\end{proof}

Let $\D\subset \Co$ be obtained by projecting any $\D(\Lambda)$. It
is independent of the choice of $\Lambda$. We consider on it the
norm given by the projection of the norm on $\D(\emptyset)$ -- any
other choice would give an equivalent norm. This is the space
described in Theorem \ref{thm:SpectralDescription}.

In a pedestrian way, the norm of an element of $\mu\in \Co$ is the
infimum of the quantity
  \begin{equation}
  \label{eq:TotalMass}
  \max\bigl(|\mu(1)|,\sup_{\p \in \Z^d} \Var(\mu_\p)\bigr)
  \end{equation}
over all decompositions $\mu=\mu(1)m^{\otimes \Z^d}+ \sum_{p\in
\Z^d} \mu_\p$ where $\mu_\p \in \B_\p$. The elements of $\D$ are
exactly those elements of $\Co$ for which such a decomposition
exists with finite \eqref{eq:TotalMass}.

\section{Proof of the main theorem}


From this point on, we fix a range $r$. The next lemma gives a
contraction estimate for the action of $P_\ve^n$ on the measures
$\mu\in\B_\p$. This is essentially contained in the paper
\cite{kl:coupling2}, one has simply to check that the only
variations involved in the computation are those of points close to
$\p$.
\begin{lem}
  \label{lem:KL}
  There exist $\nsigma,\rho\in(0,1)$, $\eps_1=\eps_1(r)>0$ and $C=C(r)>0$ such
  that, for any coupling $\Phi_\eps$ of range $r$ and strength $0\leq \eps\leq
  \eps_1$, for all $\p\in \Z^d$, for all $n\in \N$, and for all $\mu\in
  \B_\p$,
  \begin{equation}
  \label{ContractionKL}
  |P_\ve^n \mu| \leq C \nsigma^{2n}  \sup_{\jj\in \Z^d}
  \rho^{|\jj-\p|}\Var_\jj(\mu).
  \end{equation}
\end{lem}
\noindent The precise choice of the constants and the details of the
proof are provided in Section \ref{sec:proof-of-lemma}.

\begin{lem}
There exist $\nsigma\in(0,1)$ and $\eps_1=\eps_1(r)>0$
  such that, for any coupling $\Phi_\eps$ of range $r$ and strength
  $0\leq \eps\leq \eps_1$, and for any
Lipschitz function $f$ depending only on a finite number of
coordinates, there exists $C>0$ such that, for all $\mu\in \B$, for
all $n\in\N$, and for all $t\in \R$,
  \begin{equation}
  \label{LY_naive}
  \Var( P_{t,\ve}^n \mu) \leq \nsigma^{2n} \Var(\mu) + C(1+|t|) |\mu|.
  \end{equation}
\end{lem}
\begin{proof}
By Lemma \ref{lem:LipschitzControl}, we have
  \begin{equation}\label{eq:Var1}
  \Var(e^{itf}\mu) \leq \Var(\mu) + C|t| |\mu|.
  \end{equation}
Moreover, \cite[Proposition 4.1 for $\theta=1$]{kl:coupling1}
implies that
  \begin{equation}\label{eq:Var2}
  \Var( P_\ve \mu) \leq \nsigma^2 \Var(\mu) + C |\mu|.
  \end{equation}
  Using these two equations and a geometric series, we get the
  conclusion.\footnote{Referring to \cite[Proposition 4.1]{kl:coupling1} in
    this proof we make use of the assumption that $\inf|\tau'|>2$. However,
    $\var(P_{t,\eps}^\ell\mu)$ could be estimated, for each fixed $\ell$, just
    as for $\ell=1$ in \eqref{eq:Var1} and \eqref{eq:Var2}, and with
    $\nsigma\in(0,1)$ chosen such that $\nsigma^{-2\ell}<\inf|(\tau^\ell)'|$ one
    would recover \eqref{LY_naive}.}
\end{proof}

We fix from now on the value of $\eps_1(r)$ as the minimum of those
given in the two previous lemmas, it will satisfy the conclusion of
Theorem \ref{thm:SpectralDescription}. We denote also by $\nsigma$
the maximum of the values given in the previous lemmas. Fix now a
coupling $\Phi_\eps$ of range $r$ and strength
$0\leq\eps\leq\eps_1(r)$, as well as a Lipschitz function $f$
depending only on a finite number of coordinates, say coordinates in
a box $[-A,A]^d$. All the constants that we will construct from this
point on are allowed to depend on $f$ as well as $r$, $\tau$.

\begin{lem}\label{lem:ell}
  There exist $C_0>1$, $\ell>0$
  such that, for all $\p\in \Z^d$, for all $n\in \N$
  with $|\p| > \ell n$, for all $\mu\in \B_\p$,
  \begin{equation}
  \label{eq:StrongContraction}
  \Var( P_{t,\ve}^n \mu) \leq C_0(1+|t|)^2 \nsigma^{n} \Var(\mu).
  \end{equation}
\end{lem}
\begin{proof}
Write $n=a+b$ where $a=n/2$ or $(n-1)/2$ depending on whether $n$ is
even or odd. By \eqref{LY_naive},
  \begin{equation}\label{eq:xxx}
  \Var( P_{t,\ve}^{n} \mu) \leq \nsigma^{2a}
  \Var(P_{t,\ve}^b\mu) + C(1+|t|) |P_{t,\ve}^b\mu|.
  \end{equation}
For the first term, \eqref{LY_naive} again gives
$\Var(P_{t,\ve}^b\mu) \leq C(1+|t|) \Var(\mu)$. For the second one,
by \eqref{ContractionKL} and  Lemma \ref{lem:LipschitzControl},
  \begin{align*}
  |P_{t,\ve}^b\mu|
  &=
  |P_\ve^b( e^{it S_b f} \mu)|
  \leq
  C \nsigma^{2b} \sup_{\jj\in \Z^d} \rho^{|\jj-\p|}\Var_\jj( e^{it S_b f}
  \mu)
  \\&
  \leq
  C \nsigma^{2b} \sup_{\jj\in \Z^d} \rho^{|\jj-\p|} ( \Var\mu+
  \Lip_\jj(e^{itS_b f}) |\mu|).
  \end{align*}

Define a distance on $\Omega$ by $d(\x,\y)=\sup_{\q\in \Z^d} |
x_\q-y_\q|$. It does not define the product topology but,
nevertheless, there exists a constant $C$ such that $|f(\x)-f(\y)|
\leq C d(\x,\y)$ (since $f$ is Lipschitz and depends on a finite
number of coordinates). Moreover, we have $d(T_0\x,T_0\y)\leq C
d(\x,\y)$, since $\tau$ is Lipschitz, as well as $d(\Phi_\eps \x,
\Phi_\eps \y) \leq C d(\x,\y)$. Hence,
  \begin{equation}
  \Lip_\jj(e^{itS_b f}) \leq |t|\sum_{k=0}^{b-1} \Lip_\jj( f\circ T_\ve^k)
  \leq |t|\sum_{k=0}^{b-1} C^k \leq |t| C^b.
  \end{equation}
Moreover, if $|\jj|> rb+ A$, the function $e^{itS_b f}$ does not
depend on the coordinate $\jj$, hence $\Lip_\jj(e^{itS_b f})=0$.
Finally,
  \begin{equation}
  |P_{t,\ve}^b\mu| \leq C \nsigma^{2b} \Var\mu + C \nsigma^{2b}
  \sup_{|\jj|\leq rb+A} |t|\rho^{|\jj-\p|} C^b |\mu|.
  \end{equation}
  If $|\p|> \ell n$ for some large enough $\ell$, we have $\rho^{|\p| -rb-A}
  C^b \leq 1$, and we get $|P_{t,\ve}^b\mu| \leq C (1+|t|)\nsigma^{2b}
  \Var\mu$.  Together with \eqref{eq:xxx}, this proves the lemma.
\end{proof}

Fix $t\in\R$. Let $\Lambda_N$ be the box $[-\ell N,\ell N]^d$, with
$\ell$ from the previous Lemma. Let $N$ be large enough so that
$C_0(1+|t|)^2\nsigma^{N/2} \leq \frac{1}{4 B(\Lambda_N,N,r)}$. We
will now work in the extension $\E=\E(\Lambda_N)$, and study the
operator $Q=Q_{t,\ve,N,\Lambda_N}$.

\begin{lem}
There exists a constant $C>0$ such that, for all $\mu \in \E$,
  \begin{equation}
  \label{LY_extension}
  \norm{Q \mu} \leq \frac{1}{2} \norm{\mu} + C\norm{\mu}_w.
  \end{equation}
\end{lem}
\begin{proof}
By \eqref{ExpandsWeak}, we have
  \begin{equation}
  \norm{Q \mu} \leq 2B(\Lambda_N,N,r)\sup\left( \Var(P_{t,\ve}^N \mu_c),
  \sup_{\p \in \Z^d\moins \Lambda_N} \Var(P_{t,\ve}^N \mu_\p)\right).
  \end{equation}
Moreover, \eqref{eq:StrongContraction} shows that $\Var(P_{t,\ve}^N
\mu_\p) \leq C_0(1+|t|)^2 \nsigma^{N/2} \norm{\mu}$, while
\eqref{LY_naive} gives $\Var(P_{t,\ve}^N \mu_c) \leq
\nsigma^N\norm{\mu} + C(1+|t|) \norm{\mu}_w$. We get
  \begin{equation}
  \norm{Q \mu} \leq 2 B(\Lambda_N,N,r) \max(C_0(1+|t|)^2 \nsigma^{N/2} \norm{\mu},
  \nsigma^N\norm{\mu} + C(1+|t|)
  \norm{\mu}_w),
  \end{equation}
which yields the desired conclusion by the choice of $N$.
\end{proof}

This is a Lasota-Yorke inequality for the operator $Q$. The main
advantage of this construction is that, since the unit ball of
$(\E(\Lambda), \norm{\cdot})$ is relatively compact for the seminorm
$\norm{\cdot}_w$, we get from Lemma \ref{lem:hennion} that the
essential spectral radius of $Q$ is at most $1/2$. To show that the
spectral radius of $Q$ is less than $1$, it is therefore sufficient
to check that there is no eigenvalue of modulus $\geq 1$.

\begin{lem}
\label{lem:NoEigenvalue} Let $\lambda\in \C$ with $|\lambda|\geq 1$.
Let $\mu\in \E$ satisfy $Q\mu=\lambda\mu$. Then $\mu=0$, or
$|\lambda|=1$ and $tf$ is periodic.
\end{lem}
\begin{proof}
Let $\nu=\pi(\mu)\in \Co$. We will first check that it belongs to
$\B$. By Lemma \ref{lem:BoundedVarOK}, it is sufficient to check
that the variations of the measures $\nu_\Lambda$ are uniformly
bounded.

Let $\vf$ be a smooth function depending on a finite number of
coordinates, bounded by $1$, and fix $\q$. Fix $K\geq A$ such that
$|\q|\leq K$ and that $\vf$ depends only on the coordinates $\p$
with $|\p|\leq K$. For $n\in \N$ which is a multiple of $N$ we have,
since $\mu$ is an eigenfunction of $Q$,
  \begin{equation*}
  |\nu(\partial_\q \vf)|
  \leq | P_{t,\ve}^n \mu_c( \partial_\q \vf)|
  + \sum_{ |\p| \leq K+rn} | P_{t,\ve}^n \mu_\p (\partial_\q \vf)|.
  \end{equation*}
Indeed, if $|\p|>K+rn$, then
  \begin{equation}
  P_{t,\ve}^n \mu_\p( \partial_\q \vf)
  = \mu_\p( e^{it S_n f} (\partial_\q \vf) \circ T_\ve^n)
  = 0
  \end{equation}
since $e^{it S_n f} (\partial_\q \vf) \circ T_\ve^n$ does not depend
on $x_\p$. We get
  \begin{equation}
  |\nu(\partial_\q \vf)| \leq \Var(P_{t,\ve}^n \mu_c) +
  \sum_{|\p| \leq K+rn} \Var( P_{t,\ve}^n \mu_\p).
  \end{equation}
  Let $\ell'=\max(\ell, r)+1$ with $\ell$ as in Lemma~\ref{lem:ell}, and let
  $k(\p)$ be the integer part of $|\p|/\ell'$.  If $n\geq K$, then $k(\p) \leq
  n$ whenever $|\p| \leq K+rn$. Then, by \eqref{LY_naive},
  \begin{equation}
  \Var( P_{t,\ve}^n \mu_\p) \leq C(1+|t|) \Var( P_{t,\ve}^{k(\p)}
  \mu_\p).
  \end{equation}
We can then use \eqref{eq:StrongContraction} since $|\p|> \ell
k(\p)$. We get
  \begin{equation}
  \Var( P_{t,\ve}^n \mu_\p) \leq C(1+|t|)^3 \nsigma^{ k(\p)}
  \Var(\mu_\p).
  \end{equation}
Finally,
  \begin{equation}
  |\nu(\partial_\q \vf)| \leq C(1+|t|) \Var(\mu_c) + C(1+|t|)^3
  \sum_{|p| \leq K+rn} \nsigma^{ k(\p)}
  \Var(\mu_\p).
  \end{equation}
This last sum is bounded uniformly in $K$ and $n$. This proves that
the variation of the measures $\nu_\Lambda$ are uniformly bounded,
i.e. $\nu\in \B$.

If $\nu=0$, the marginal of $\nu$ on $\Lambda_N$ vanishes, i.e.,
$\mu_c=0$. Therefore, $\norm{\mu}_w=0$. The Lasota-Yorke inequality
\eqref{LY_extension} then gives $\mu=0$.

Assume now that $\nu\not=0$. We will prove that $|\lambda|=1$ and
that $tf$ is periodic. The measure $\nu$ satisfies $P_{t,\ve}^N\nu=
\lambda\nu$. The absolute value $A(\nu)$ of the measure $\nu$
belongs to $\B$, by Lemma \ref{lem:AinB}. It satisfies
  \begin{equation}
  A(\nu) = |\lambda|^{-1} A( P_{t,\ve}^N \nu) \leq P_\ve^N A(\nu),
  \end{equation}
where the last inequality is obtained by the following direct
computation:
\begin{equation}
  \begin{split}
    A( P_{t,\ve}^N \nu)(\varphi)
    &=\sup_{|g|\leq1}|( P_{t,\ve}^N \nu)(g\cdot\varphi)|
    =\sup_{|g|\leq1}\left|\nu\big(e^{itS_Nf}\cdot g\circ T_\eps^N\cdot\varphi\circ T_\ve^N
      \big)\right|
    \\&\leq \sup_{|g|\leq1}|\nu(g\cdot \varphi\circ T_\ve^N )|
    =P_\ve^N A(\nu)(\varphi).
  \end{split}
\end{equation}

Since $A(\nu)$ and $P_\ve^N A(\nu)$ have the same mass, this yields
$A(\nu)=P_\ve^N A(\nu)$ and $|\lambda|=1$. Since $A(\nu)$ belongs to
$\B$, it has to be a scalar multiple of the SRB measure $\mu_\ve$,
see Theorem~\ref{Thm:InvMeas}. In particular, $\nu$ is absolutely
continuous with respect to $\mu_\ve$, and the Radon-Nikodym
derivative $g= \frac{ \dd \nu}{\dd
  \mu_\ve}$ is a function of almost everywhere constant modulus.  Since we
assume $\nu$ to be nonzero, we have $|g|\neq0$ almost everywhere.
Then
  \begin{equation}
  P_\ve^N \left( e^{itS_N f} \frac{g}{g\circ T_\ve^N} \mu_\ve\right)
  = \frac{1}{g} P_\ve^N(e^{itS_N f} g \mu_\ve)
  = \frac{1}{g} P_{t,\ve}^N( \nu)
  = \lambda \frac{1}{g} \nu
  = \lambda \mu_\ve.
  \end{equation}
In particular,
  \begin{equation}
  1= \left| \int e^{itS_N f} \frac{g}{g\circ T_\ve^N} \dd\mu_\ve\right|
  \leq \int \left| e^{itS_N f} \frac{g}{g\circ T_\ve^N}\right|\dd\mu_\ve = 1.
  \end{equation}
Therefore, we have equality in the inequality, and $e^{itS_N f}
\frac{g}{g\circ T_\ve^N}$ is almost everywhere equal to a constant
of modulus $1$. This shows that $tS_N f$ is periodic for $T_\ve^N$.
By Lemma \ref{lem:StableCohom}, $tf$ is periodic for $T_\ve$.
\end{proof}

\begin{proof}[Proof of Theorem \ref{thm:SpectralDescription}]

The extension $\D$ is as described in Paragraph
\ref{subsec:DescribesD}. The formula \eqref{eq:TotalMass} for the
norm clearly gives $|\mu(1)|\leq \norm{\mu}$.

Let $\Lambda$ be a fixed box, we want to check equation
\eqref{eq:Algebra}, i.e.
  \begin{equation}
  \label{eq:ContMult}
  \norm{u\mu}\leq C(\Lambda) (\Lip(u)+|u|_\infty)\norm{\mu}
  \end{equation}
whenever $u$ is Lipschitz continuous and depends only on coordinates
in $\Lambda$. Let us first work in the extension $\E(\Lambda)$. If
$\mu=(\mu_c, (\mu_\p)_{p\in \Z^d\moins \Lambda}) \in \E(\Lambda)$,
then all the measures $u\mu_\p$ still belong to $\B_\p$, and they
satisfy $\Var(u\mu_\p) \leq (\Lip u+|u|_\infty) \Var \mu_\p$ by
Lemma \ref{lem:LipschitzControl}. Moreover, $u\mu_c$ is still of the
form $\nu \otimes m^{\Z^d\moins \Lambda}$, and its variation is at
most $(\Lip u+|u|_\infty) \Var \mu_c$, again by Lemma
\ref{lem:LipschitzControl}. Hence, the multiplication by $u$ is well
defined on $\E(\Lambda)$ and its norm is at most $(\Lip
u+|u|_\infty)$. This multiplication leaves $\E(\Lambda)_0$
invariant, hence induces an operator on the quotient space
$\D(\Lambda)$ with the same bound on its norm (see
section~\ref{subsec:DescribesD}). Since $\D(\Lambda)$ is isomorphic
to $\D$ by Proposition \ref{prop:EquivNorms}, this proves
\eqref{eq:ContMult}.

In order to get analyticity of $t\mapsto P_{t,\ve}$, it is enough to
prove that the map $M_t(\mu)=e^{itf}\mu$ depends analytically on
$t$. For this, we only have to check that the series expansion
  \begin{equation}
  \sum_{n\geq0}\frac{(itf)^n}{n!}\mu
  \end{equation}
is well defined for any $\mu\in\D$. But this is a direct consequence
of \eqref{eq:ContMult}, since
  \begin{equation}
  \begin{split}
  \norm{\frac{(itf)^n}{n!}\mu}&\leq C(\Lip(f^n)+|f^n|_\infty)\norm{\mu}
  \frac{|t|^n}{n!}
  \\&\leq C(n\Lip(f)+|f|_\infty)|f|_\infty^{n-1}\norm{\mu}\frac{|t|^n}{n!}.
  \end{split}
  \end{equation}
This gives analyticity of $M_t$, and its series expansion.

In \cite{kl:coupling2}, it is proved that, in the extension
$\E(\emptyset)$, the operator $Q_{0,\ve,N,\emptyset}$ (which is a
lift of $P_{0,\ve}^N$ on $\Co$) has a simple eigenvalue at $1$ for
sufficiently large $N$, the rest of its spectrum being contained in
a disk of radius $<1$: this is indeed an easy consequence of Lemma
\ref{lem:KL}. After a quotient by $\E(\emptyset)_0$ (which is left
invariant by $Q_{0,\ve,N,\emptyset}$), this implies that $P_\ve^N$
acts continuously on $\D$, has a simple eigenvalue at $1$ and the
rest of its spectrum is contained in a disk of radius $<1$. The same
is then true for the operator $P_\eps$ itself.

Consider now $t\not=0$ and assume that $f$ is aperiodic. For a
suitable $N$ and a suitable box $\Lambda_N$, Lemma
\ref{lem:NoEigenvalue} shows that the spectral radius of
$Q_{t,\ve,N,\Lambda_N}$ is $<1$ on $\E(\Lambda_N)$. On the quotient
$\D\cong \E(\Lambda_N)/ \E(\Lambda_N)_0$, this implies that the
spectral radius of $P_{t,\ve}^N$ is $<1$. Therefore, $P_{t,\ve}$
also has a spectral radius $<1$.
\end{proof}

\section{Proof of Lemma~\ref{lem:KL}}
\label{sec:proof-of-lemma}

We introduce a family of additional ``local'' norms: for
$\rho\in(0,1)$ (to be chosen later), for any $\p\in\Z^d$,
$\Lambda\subset\Z^d$ and $\mu\in\B$ let
\begin{align}
  \label{eq:local-var}
  \Var^\p(\mu)&=\sup_{\jj\in\Z^d}\rho^{|\jj-\p|}\Var_\jj(\mu)\;,\\
  \Var^\Lambda(\mu)&=\sup_{\p\in\Lambda}\Var^\p(\mu)\;.
\end{align}
Observe that $\Var(\mu)=\sup_\p\Var^\p(\mu)=\Var^{\Z^d}(\mu)$.  In
this section we denote $\LL{\p}{n}=\{\q \st |\q-\p|\leq rn\}$, so
the range $r$ will often be suppressed in the notation.

For the proof of Lemma~\ref{lem:KL} we need two further lemmas that
will be proved later. Let $\lambda_1=\frac12\inf|\tau'|>1$ and
denote by $\nsigma_0\in(0,1)$ the mixing rate of $\tau$.

\begin{lem}[Localized Lasota-Yorke type estimate]
  \label{lem:LY-local}
  For any $\lambda\in(1,\lambda_1)$, for any range $r$ and any $\rho\in(0,1)$,
  there are $\ve_2>0$ and $C>0$ such that, for any coupling $\Phi_\eps$ of
  range $r$ and strength $0\leq \eps\leq \ve_2$, for all $m\in\N$, for all
  $\p\in\Z^d$, and for all $\nu\in\B$,
  \begin{equation}
    \label{eq:gkLY}
    \Var^\p(P_\ve^m\nu)
    \leq
    C\left(\lambda^{-m}\Var^\p(\nu)+|\nu|\right)
    \leq
    C\Var^\p(\nu)\;.
  \end{equation}
\end{lem}
\begin{lem}\label{lem:eps}
For any range $r$ and any $\rho\in(0,1)$, there are $\ve_3>0$ and
$C>0$ such that, for any coupling $\Phi_\eps$ of range $r$ and
strength $0\leq \eps\leq \ve_3$, for all $m\in\N$, for all
$\p\in\Z^d$, and for all $\nu\in\B_\p$,
  \begin{equation}
    \label{eq:gkeps}
    |P_\ve^m\nu|
    \leq
    C\,\left(\nsigma_0^m\Var_\p(\nu)+m\ve\Var^\p(\nu)\right)\;.
  \end{equation}
\end{lem}

\begin{proof}[Proof of Lemma~\ref{lem:KL}]
We can precise the choice of the constants appearing in the lemma:
let $\lambda\in(1,\lambda_1)$ be fixed, then choose
$\nsigma,\nsigma_1,\rho\in(0,1)$ such that
\begin{equation}
  \label{eq:rates}
  \sqrt{\max\{\lambda^{-1},\nsigma_0\}}<\nsigma_1<\rho^{2r}\nsigma\;.
\end{equation}
The maximal coupling strength $\ve_1$ will have to be taken smaller
than $\ve_2$ and $\ve_3$ from the previous lemmas, and even smaller
in the calculation below.

Before getting into the proof of Lemma \ref{lem:KL}, let us
establish a preliminary inequality in the extension $\E(\emptyset)$,
using  Lemmas \ref{lem:LY-local} and \ref{lem:eps}. Let
$Q:=Q_{t,\ve,2m,\emptyset}$ be the lift of $P_{t,\eps}^{2m}$
described in section~\ref{subsec:BV}. It redistributes mass from a
site $\q$ to sites in $\LL{\q}{2m}$ only. We claim that there exist
$m\in \N$ and $\ve_1>0$ such that, whenever the coupling strength is
at most $\ve_1$, for each $\Gamma\subseteq\Z^d$ and each
$\tilde\nu\in\E(\emptyset)$ with $\tilde\nu_c=0$,
\begin{equation}
  \label{eq:gk2}
  \sup_{\jj\in\q+\Gamma}\Var^{\q+\Gamma}\big((Q{\tilde\nu})_\jj\big)
  \leq
  \nsigma_1^{2m}\,\sup_{\jj\in\LL{\q}{2m}+\Gamma}
  \Var^{\LL{\q}{2m}+\Gamma}({\tilde\nu}_\jj)
\end{equation}

In view of the redistribution mechanism described in
\eqref{ExpandsWeak}, we have
\begin{equation}
  \label{eq:gk3}
  \begin{split}
    \sup_{\jj\in\q+\Gamma}\Var^{\q+\Gamma}\big(&(Q{\tilde\nu})_\jj\big)\\
    \leq&
    C\,m^d\,\sup_{\jj\in\LL{\q}{2m}+\Gamma}\Var^{\q+\Gamma}(P_\ve^{2m}{\tilde\nu}_\jj)\\
    \leq&
    C\,m^d\,\sup_{\ii,\jj\in\LL{\q}{2m}+\Gamma}\left(
      \lambda^{-m}\Var^\ii(P_\ve^{m}{\tilde\nu}_\jj)+
      |P_\ve^m{\tilde\nu}_\jj|\right)\\
    \leq&
    C\,m^d\,\sup_{\ii,\jj\in\LL{\q}{2m}+\Gamma}\big(
      \lambda^{-m}\Var^\ii({\tilde\nu}_\jj)\\
      &\hspace*{30mm}+\nsigma_0^m\Var_\jj({\tilde\nu}_\jj)+m\ve\,\Var^\jj({\tilde\nu}_\jj)
      \big)
    \end{split}
\end{equation}
where we used the Lasota-Yorke type inequality \eqref{eq:gkLY} and
the estimate \eqref{eq:gkeps}. Hence,
\begin{equation}
  \label{eq:gk4}
  \begin{split}
    &\sup_{\jj\in\q+\Gamma}\Var^{\q+\Gamma}\big((Q{\tilde\nu})_\jj\big)\\
    &\hspace*{1cm}\leq
    C\,m^d\,\sup_{\ii,\jj\in\LL{\q}{2m}+\Gamma}
    \big(\lambda^{-m}\Var^\ii(\tilde\nu_\jj)
    +(\nsigma_0^m+m\ve)\Var^\jj(\tilde\nu^\jj)\big)\\
    &\hspace*{1cm}\leq
    C\,m^d\,\big(\lambda^{-m}+\nsigma_0^m+m\ve\big)
    \sup_{\jj\in\LL{\q}{2m}+\Gamma}\Var^{\LL{\q}{2m}+\Gamma}(\tilde\nu_\jj).
  \end{split}
\end{equation}
Choosing $m$ sufficiently large and then $\ve_1$ sufficiently small,
\eqref{eq:gk2} follows.

Let us now prove Lemma \ref{lem:KL}. As $P_\ve$ contracts the total
variation norm, it suffices to prove the lemma for multiples $n=k2m$
of the fixed integer $2m$ which satisfies \eqref{eq:gk2}. So let
$\mu\in\B_\p$, define an element $\tilde\mu$ which has only zero
components except for $\tilde\mu_\p=\mu$ and observe first that
\begin{equation}
  \label{eq:gk1}
  \begin{split}
    |P_\ve^{k2m}\mu|
    &\leq
    \sum_{\q\in\Z^d}\Big|\big(Q^k \tilde\mu\big)_\q\Big|
    \leq
    \frac12\sum_{\q\in\LL{\p}{k2m}}
    \Var_\q\left(\big(Q^k\tilde\mu\big)_\q\right)\\
     &\leq
     C\cdot(k2m)^d\cdot\sup_{\q\in\LL{\p}{k2m}}
     \Var_\q\left(\big(Q^k\tilde\mu\big)_\q\right)
    \;,
  \end{split}
\end{equation}
where we used that each application of $Q:=Q_{t,\ve,2m,\emptyset}$
redistributes mass from a site $\q$ to sites in $\LL{\q}{2m}$ only.

Applying \eqref{eq:gk2} repeatedly and observing that
$\tilde\mu_\jj=0$ if $\jj\neq\p$ and $\tilde\mu_\p=\mu$, we obtain
\begin{equation}
  \label{eq:gk5}
  \begin{split}
    \sup_{\q\in\LL{\p}{k2m}}\Var_{\q}\big((Q^k{\tilde\mu})_\q\big)
    &\leq
    \sup_{\q\in\LL{\p}{k2m}}\nsigma_1^{k2m}\sup_{\jj\in\LL{\q}{k2m}}\Var^{\LL{\q}{k2m}}(\tilde\mu_\jj)\\
    &=
    \nsigma_1^{k2m}\,\sup_{\q\in\LL{\p}{k2m}}
    \Var^{\LL{\q}{k2m}}({\tilde\mu}_\p)\\
    &\leq
    \nsigma_1^{k2m}\rho^{-2rk2m}\Var^\p(\mu)\;.
  \end{split}
\end{equation}
Together with \eqref{eq:gk1} this yields $|P_\ve^n\mu|\leq C\, n^d
(\nsigma_1\rho^{-2r})^n \Var^\p(\mu)$, which finishes the proof of
Lemma~\ref{lem:KL} in view of the choice of the constants in
\eqref{eq:rates}.
\end{proof}

\begin{proof}[Proof of Lemma~\ref{lem:LY-local}]
  We will prove
  \begin{equation}
    \label{eq:LYP1}
    \Var^\p(P_\ve\nu)
    \leq
    \lambda^{-1}\Var^\p(\nu)+C|\nu|\;.
  \end{equation}
  From this \eqref{eq:gkLY} follows by induction.\footnote{It is only
    here where we use the assumption $\inf|\tau'|>2$. For
    $1<\inf|\tau'|\leq2$ this reduction to the case $m=1$ is not
    possible, see also the remarks in \cite[Footnote
    12]{kl:coupling2}.}

  Observe first that
  \begin{equation}
    \label{eq:LYP2}
    \Var^\p(P_0\nu)
    \leq
    \lambda_1^{-1}\Var^\p(\nu)+C|\nu|\;,
  \end{equation}
  where $\lambda_1=\frac12\inf|\tau'|$. This is a simple consequence
  of the Lasota-Yorke inequality for the single site map, compare
  e.g.~the proof of Lemma 3.2 in \cite{kl:coupling1}. We will show
that\footnote{We write $F_*\nu$ for the push-forward of the measure
$\nu$ under the map $F$, i.e., the measure given by
$(F_*\nu)(A)=\nu(F^{-1}A)$. Note that this object is sometimes
denoted by $F^*\nu$ in \cite{kl:coupling1, kl:coupling2}.}
  \begin{equation}
    \label{eq:LYP3}
    \Var^\p((\Phi_\ve)_*\nu)
    \leq
    (1+C\ve)\Var^\p(\nu)\;.
  \end{equation}
We first notice that under a mild bound on the coupling strength
$\ve$,  the coupling assumption
  \eqref{eq:couplingass} ensures that the infinite matrix $D\Phi_\ve(\x)$ is
  invertible. Moreover, taking $C$ large enough to get
$1_{|\ii-\jj|\leq r}\leq \frac C2\rho^{4|\ii-\jj|}$ for all
$\ii,\jj\in\Z^d$, the second and third part of this assumption can
be rewritten as
\begin{equation}
 |(DA_\eps)_{\ii\jj}|_\infty\le C\eps\rho^{4|\ii-\jj|},\quad
  |\partial_{\kk}(DA_\eps)_{\ii\jj}|_\infty\le C\eps\rho^{4|\ii-\jj|}\;.
\end{equation}
A direct computation using these estimates (see for example page 300
in \cite{jiapes}) gives that $B(\x):=(D\Phi_\ve(\x))^{-1}$ satisfies
  \begin{equation}
    |b_{\ii\ii}|_\infty\leq1+C\ve,\quad
    |b_{\ii\jj}|_\infty\leq C\ve\rho^{2|\jj-\ii|},\quad
    |\partial_\ii b_{\ii\jj}|_\infty\leq C\ve\rho^{2|\jj-\ii|}\;.
  \end{equation}
We can then follow the proof of Lemma 3.3 in \cite{kl:coupling1}
with
  some modifications. For all $\jj,\p\in\Z^d$,
  \begin{equation*}
    \label{eq:LYP4}
    \begin{split}
      \rho^{|\jj-\p|}&((\Phi_\ve)_*\nu)(\partial_\jj\vf)\\
      &=
      \rho^{|\jj-\p|}\sum_{\ii\in\Z^d}
      \nu\big(\partial_\ii(\vf\circ\Phi_\ve)\,b_{\ii\jj}\big)\\
      &=
      \rho^{|\jj-\p|}\sum_{\ii\in\Z^d}
      \nu\big(\partial_\ii(\vf\circ\Phi_\ve\cdot
      b_{\ii\jj})\big)
      -\rho^{|\jj-\p|}\sum_{\ii\in\Z^d}
      \nu\big(\vf\circ\Phi_\ve\cdot\partial_\ii b_{\ii\jj}\big)\\
      &\leq
      \rho^{|\jj-\p|}\sum_{\ii\in\Z^d}
      \Var_\ii(\nu)\cdot|b_{\ii\jj}|_\infty
      +|\nu|\rho^{|\jj-\p|}\cdot\sum_{\ii\in\Z^d}
      |\partial_\ii b_{\ii\jj}|_\infty\\
      &\leq
      \rho^{|\jj-\p|}\Var_\jj(\nu)+
        C\,\ve\Big(\Var^\p(\nu)\sum_{\ii\in\Z^d}
        \rho^{|\jj-\p|-|\ii-\p|+2|\jj-\ii|}
        +|\nu|\sum_{\ii\in\Z^d}\rho^{|\jj-\p|+2|\jj-\ii|}\Big)\\
        &\leq
        (1+C\ve)\Var^\p(\nu)\,.
    \end{split}
  \end{equation*}
  This yields \eqref{eq:LYP3} and finishes the proof of
  Lemma~\ref{lem:LY-local}.
\end{proof}

\begin{proof}[Proof of Lemma~\ref{lem:eps}]
  The proof follows closely the corresponding one in \cite{kl:coupling2}. For
  each $\p\in\Z^d$ define a coupling map $\Phi_{\ve,\p}:\Omega\to\Omega$ where
  site $\p$ is decoupled from all other sites,
  \begin{equation}
    \label{eq:gk8}
    (\Phi_{\ve,\p}(\x))_\q
    =
    \begin{cases}
      x_\p&\text{if }\q=\p\\
      (\Phi_{\ve}(\x_{\Z^d\setminus\{\p\}},a))_\q&\text{if }\q\neq\p
    \end{cases}
  \end{equation}
  where $a$ is an arbitrary point in $I$. Denote by $P_{\ve,\p}$ the
  Perron-Frobenius operator of $\Phi_{\ve,\p}\circ T_0$. We will show
  that, for each
  $\nu\in\B$,
  \begin{equation}
    \label{eq:gk9}
    |(\Phi_\ve)_*\nu-(\Phi_{\ve,\p})_*\nu|
    \leq
    C\ve\Var^\p(\nu)\;.
  \end{equation}
  Then, making use of the fact that $|P_\ve|=|P_{\ve,\p}|=1$ and of estimate
  \eqref{eq:gkLY}, a simple telescoping argument yields
  \begin{equation}
    \label{eq:gk10}
    |P_\ve^m\nu-P_{\ve,\p}^m\nu|
    \leq
    C\,m\ve\Var^\p(\nu)\;,
  \end{equation}
  and \eqref{eq:gkeps} follows once we have shown that
  $|P_{\ve,\p}^m\nu|\leq C\nsigma_0^m\Var_\p(\nu)$ for any
  $\nu\in\B_\p$.
  But this is proved
  precisely as in \cite[p.40/41]{kl:coupling2}, where $\nsigma_0$ is
  the mixing rate for the single site map.

  It remains to prove \eqref{eq:gk9}. Here we can follow closely the
  proof of Lemma 3.2a) in \cite{kl:coupling2}. Indeed, let
  $F_t:=t\Phi_{\ve,\p}+(1-t)\Phi_\ve$ and
  $\Delta_\q:=(\Phi_{\ve,\p}-\Phi_\ve)_\q$. Just as in
  \cite{kl:coupling2} one shows that, for each test function $\vf$,
  \begin{equation}
    \label{eq:gk11}
    \big((\Phi_{\ve,\p})_*\nu-(\Phi_\ve)_*\nu\big)(\vf)
    =
    \int_0^1\sum_{\q\in\Z^d}
    (F_t)_*\big(\Delta_\q\cdot\nu\big)(\partial_\q\vf)\,dt\;.
  \end{equation}
  As, in our case, $\Delta_\q=0$ if $|\q-\p|> r$, we conclude
  \begin{equation}
    \label{eq:gk12}
    \begin{split}
      |(\Phi_{\ve,\p})_*\nu-(\Phi_\ve)_*\nu|
      &\leq
      \sum_{|\q-\p|\leq r}\sup_{0\leq t\leq1}
      \Var_\q\Big((F_t)_*\big(\Delta_\q\cdot\nu\big)\Big)\\
      &\leq
      C\sum_{|\q-\p|\leq r}\Var^\q\Big(\Delta_\q\cdot\nu\Big)\;,
    \end{split}
  \end{equation}
  where we used \eqref{eq:LYP3} (which applies as well to $(F_t)_*$) for
  the second inequality. Hence, by Lemma \ref{lem:LipschitzControl},
  \begin{align*}
  |(\Phi_{\ve,\p})_*\nu-(\Phi_\ve)_*\nu|
  &
  \leq C\sum_{|\q-\p|\leq r}\left(
  |\Delta_\q|_\infty\Var^\q(\nu)+
  \sup_{\jj\in\Z^d}\rho^{|\jj-\q|}\Lip_j(\Delta_\q)|\nu|\right)
  \\ & \leq C\ve \Var^\q(\nu)
  \end{align*}
  in view of assumption
  \eqref{eq:couplingass}. This is \eqref{eq:gk9} and finishes the proof of the
  lemma.
\end{proof}

\appendix

\section{Proof of Theorem \ref{MainThm} assuming Theorem
\ref{thm:SpectralDescription}} 
\label{append}

The operator $P_\ve$ has a simple
eigenvalue at $1$, and the corresponding eigenfunction is the
invariant measure $\mu_\ve$ obtained in \cite{kl:coupling2}. By
classical analytic perturbation theory, the operator $P_{t,\ve}$ has
for small $t$ a unique eigenvalue $\lambda(t)$ close to $1$, which
is still simple. Let $\Pi_t$ denote the corresponding spectral
projection, and $\mu_{t,\ve}= \Pi_t(\mu_\ve)$. There exist
$\delta<1$ and $C>0$ such that, for all small enough $t$, for all
$n\in \N$,
  \begin{equation}
  \label{SpectralDescription}
  \left| \int e^{it S_n f} \dd\mu_\ve - \lambda(t)^n \mu_{t,\ve}(1)
  \right| \leq C \delta^n.
  \end{equation}
Hence, a precise description of the eigenvalue $\lambda(t)$ will
imply a central limit theorem for the Birkhoff sums $S_n f$.

Let
$\nu_t=\mu_{t,\ve}/ \mu_{t,\ve}(1)$. Differentiating the equality
$P_{t,\ve} \nu_t= \lambda(t) \nu_t$ and using $\left.\frac{d
P_{t,\ve}}{dt}\right|_{t=0}= P_{\ve}(if \cdot)$, we get
  \begin{equation}
  \label{Diff1}
  P_\ve(if\mu_\ve)+ P_\ve( \nu'_0) = \lambda'(0) \mu_\ve + \nu'_0.
  \end{equation}
Integrating the function $1$ with respect to this equality, we
obtain
  \begin{equation}
  i \int f\dd\mu_\ve + \int \dd \nu'_0 = \lambda'(0) + \int \dd\nu'_0.
  \end{equation}
Since $\int f\dd\mu_\ve=0$, we therefore have $\lambda'(0)=0$.

Differentiating twice $P_{t,\ve} \nu_t= \lambda(t) \nu_t$ yields
  \begin{equation}
  P_\ve(-f^2 \mu_\ve)+2 P_\ve( if \nu'_0) + P_\ve( \nu''_0)
  = \nu''_0 + \lambda''(0) \mu_\ve.
  \end{equation}
Integrating the function $1$ yields
  \begin{equation}
  \label{Diff2}
  \lambda''(0)= -\int f^2 \dd\mu_\ve + 2i \int f \dd\nu'_0.
  \end{equation}
>From \eqref{Diff1}, we have $\nu'_0= P_\ve \nu'_0 + P_\ve(if
\mu_\ve)$. Iterating this equation gives
  \begin{equation}
  \label{limit}
  \nu'_0= P_\ve^n \nu'_0 + i \sum_{k=1}^n P_\ve^k (f \mu_\ve).
  \end{equation}
Since $\nu_t(1)=1$, we have $\nu'_0(1)=0$. The space $\{ \mu\in
\D\st \mu(1)=0\}$ is closed and $P_\ve$ leaves this space invariant,
therefore its spectral radius on this space is $<1$. This implies
that $P_\ve^n \nu'_0$ converges exponentially fast to $0$. In the
same way, $(f\mu_\ve)(1)=0$, hence $P_\ve^k(f \mu_\ve)$ converges
exponentially fast to $0$ in $\D$. Letting $n$ tend to infinity in
\eqref{limit}, we get
$  \nu'_0= i\sum_{k=1}^\infty P_\ve^k( f\mu_\ve)$.
In particular,
$  \nu'_0(f)= i \sum_{k=1}^\infty \int f\cdot f\circ T_\ve^k
\dd\mu_\ve$,
and this series converges exponentially fast. From \eqref{Diff2}, we
obtain
  \begin{equation}
   \lambda''(0) = - \int f^2 \dd\mu_\ve - 2 \sum_{k=1}^\infty \int
  f\cdot f\circ T_\ve^k \dd\mu_\ve.
  \end{equation}
Moreover,
  \begin{align*}
  \int \left(\sum_{k=0}^{n-1} f\circ T_\ve^k \right)^2 \dd\mu_\ve
  &
  = n\int f^2\dd\mu_\ve
  + 2\sum_{k=1}^{n} (n-k) \int f\cdot f\circ T_\ve^k \dd\mu_\ve
  \\&
  = -n \lambda''(0) - 2 \sum_{k=1}^\infty k \int f\cdot f\circ T_\ve^k
  \dd\mu_\ve + O(\delta^n)
  \\&
  = - n\lambda''(0) + O(1).
  \end{align*}
Since this integral is nonnegative, this shows that
$\lambda''(0)\leq 0$. Hence, we can write $\lambda''(0)=-\sigma^2$
for some $\sigma\geq 0$. Furthermore, if $\lambda''(0)=0$, then $S_n
f$ is bounded in $L^2$, which implies that $f$ can be written as
$u-u\circ T_\ve$ in $L^2$ (see e.g. \cite{vitesse_birkhoff}). This
proves the non-degeneracy criterion in Theorem \ref{MainThm}.

Since $\lambda(t)=- \sigma^2 t^2/2 + o(t^2)$, $\lambda(t/\sqrt{n})^n
\to e^{-\sigma^2 t^2/2}$. Together with \eqref{SpectralDescription},
this shows that $S_n f/\sqrt{n}$ converges in distribution to
$\boN(0,\sigma^2)$ and proves the central limit theorem. The local
limit theorem is then easily derived from the description of
$\lambda(t)$ for small $t$ and the control of the spectral radius of
$P_{t,\ve}$ for all $t\not=0$, see \cite{hennion_herve} for further
details.

\bibliography{biblio}
\bibliographystyle{alpha}

\end{document}